\documentclass[a4paper,11pt]{amsart}
\usepackage{amsfonts,amsmath,amsthm,amssymb,booktabs}
\usepackage{enumitem}
\usepackage{mathrsfs}
\usepackage[margin=1in]{geometry}
\usepackage[colorlinks=true,linkcolor=blue,filecolor=blue,citecolor=red]{hyperref}
\usepackage{orcidlink}
\usepackage{cleveref}

\theoremstyle{plain}
\newtheorem{theorem}{Theorem}[section]
\newtheorem{corollary}[theorem]{Corollary}
\newtheorem{lemma}[theorem]{Lemma}

\newtheorem{proposition}[theorem]{Proposition}

\theoremstyle{definition}

\newtheorem{example}[theorem]{Example}

\begin{document}
\title[Certain products of the Rogers-Ramanujan continued fraction]{Modularity of certain products of the Rogers-Ramanujan continued fraction}
\author[Russelle Guadalupe]{Russelle Guadalupe\orcidlink{0009-0001-8974-4502}}
\address{Institute of Mathematics, University of the Philippines Diliman\\
Quezon City 1101, Philippines}
\email{rguadalupe@math.upd.edu.ph}

\renewcommand{\thefootnote}{}



\renewcommand{\thefootnote}{\arabic{footnote}}
\setcounter{footnote}{0}

\begin{abstract}
We study the modularity of the functions of the form $r(\tau)^ar(2\tau)^b$, where $a$ and $b$ are integers with $(a,b)\neq (0,0)$ and $r(\tau)$ is the Rogers-Ramanujan continued fraction, which may be considered as companions to the Ramanujan's function $k(\tau)=r(\tau)r(2\tau)^2$. In particular, we show that under some condition on $a$ and $b$, there are finitely many such functions generating the field of all modular functions on the congruence subgroup $\Gamma_1(10)$. Furthermore, we establish certain arithmetic properties of the function $l(\tau)=r(2\tau)/r(\tau)^2$, which can be used to evaluate these products. We employ the methods of Lee and Park, and some properties of $\eta$-quotients and generalized $\eta$-quotients to prove our results.
\end{abstract}

\maketitle

\section{Introduction}\label{sec1}

For an element $\tau$ in the complex upper half-plane $\mathbb{H}$ and $q := e^{2\pi i\tau}$, the Rogers-Ramanujan continued fraction is defined by
\begin{align*}
r(\tau) := \dfrac{q^{1/5}}{1+\dfrac{q}{1+\dfrac{q^2}{1+\dfrac{q^3}{1+\cdots}}}} = q^{1/5}\prod_{n=1}^{\infty}\dfrac{(1-q^{5n-1})(1-q^{5n-4})}{(1-q^{5n-2})(1-q^{5n-3})}.
\end{align*}

In his first correspondence \cite[pp. 21--30]{ramlet} to Hardy in 1913, Ramanujan asserted without proof that
\begin{align*}
r(i) = \sqrt{\dfrac{5+\sqrt{5}}{2}}-\dfrac{\sqrt{5}+1}{2}.
\end{align*}
He later claimed in his second letter \cite[pp. 53--62]{ramlet} to Hardy that the values of $r(\frac{\sqrt{-n}}{2})$ \textquotedblleft can be exactly found\textquotedblright\,for any positive rational number $n$. Ramanujan's observations were settled by Gee and Honsbeek \cite{geehons}, who treated $r(\tau)$ as a modular function, and used Watson's formulas \cite{wat} to show that $r(\tau)$ generates the field of all modular functions on the congruence subgroup $\Gamma(5)$. Consequently, they determined the exact values of $r(\tau)$ at some imaginary quadratic points $\tau\in\mathbb{H}$ as nested radicals.

Recently, Lee and Park \cite{leepark} proved that the Ramanujan's function $k(\tau):=r(\tau)r(2\tau)^2$ generates the field of all modular functions on the congruence subgroup $\Gamma_1(10)$. They also proved that if $K$ is an imaginary quadratic field with ring of integers $\mathbb{Z}[\tau]$ and $\tau\in K\cap\mathbb{H}$, then $k(\tau)$ generates the ray class field $K_{(10)}$ modulo $10$ over $K$. As a result, they computed the exact values of $r(\tau)$ at some imaginary quadratic points $\tau\in\mathbb{H}$. We refer the reader to \cite{andber,coop,xiayao,ye} for other arithmetic properties of $k(\tau)$.

In this paper, we partially extend Lee and Park's result by exploring the functions $r(\tau)^ar(2\tau)^b$ for integers $a$ and $b$ with $(a,b)\neq (0,0)$, which may be considered as companions of $k(\tau)$. More precisely, we show that under some condition on $a$ and $b$, there are finitely many such functions generating the field of all modular functions on $\Gamma_1(10)$. We also show that under the same condition, the above functions can be expressed in terms of $l(\tau) := r(2\tau)/r(\tau)^2$. We further establish certain arithmetic properties of $l(\tau)$, which will be useful in computing the exact values of $r(\tau)^ar(2\tau)^b$ at imaginary quadratic points $\tau\in\mathbb{H}$. We apply the methods of Lee and Park \cite{leepark}, as well as some properties of $\eta$-quotients and generalized $\eta$-quotients defined by Yang \cite{yang}, to prove our results.

We organize the paper as follows. In Section \ref{sec2}, we give a brief background about modular functions, particularly $\eta$-quotients and generalized $\eta$-quotients, on some congruence subgroups of $\mbox{SL}_2(\mathbb{Z})$. In Section \ref{sec3}, we discuss the modularity of the functions $r(\tau)^ar(2\tau)^b$ on $\Gamma_1(10)$, where $a$ and $b$ are integers with $(a,b)\neq (0,0)$ and $a\equiv 3b\pmod{5}$. We show that these functions can be expressed as a rational function of $l(\tau)$. In Section \ref{sec4}, we prove that modular equations for $l(\tau)$ of any level always exist, and derive such equations for $l(\tau)$ of levels two and five.  We prove in Section \ref{sec5} that if $K$ is an imaginary quadratic field with ring of integers $\mathbb{Z}[\tau]$ for some $\tau\in K\cap\mathbb{H}$, then $l(\tau)$ generates $K_{(10)}$. We also prove that $l(\tau)$ and the functions $r(\tau)^ar(2\tau)^b$ are algebraic units for $\tau\in K\cap\mathbb{H}$. Finally, we give an example of finding the class polynomial for $K_{(10)}$ when the discriminant of $K$ is divisible by $4$ using the Shimura reciprocity law due to Cho and Koo \cite{chokoo}, Gee \cite{gee}, and the author \cite{guad} to evaluate $r(\tau)^ar(2\tau)^b$.  We have performed most of our calculations via \textit{Mathematica}.

\section{Preliminaries}\label{sec2}

In this section, we recall some facts about modular functions on some congruence subgroups, which are subgroups of $\mbox{SL}_2(\mathbb{Z})$ containing 
\begin{equation*}
\Gamma(N) := \left\lbrace\begin{bmatrix}
a & b\\ c& d
\end{bmatrix} \in \mbox{SL}_2(\mathbb{Z}) : a-1\equiv b\equiv c\equiv d-1\equiv 0\pmod N\right\rbrace
\end{equation*}
for some integer $N\geq 1$. Examples of congruence subgroups we focus on in this paper are
\[\begin{aligned}
\Gamma_1(N) &:= \left\lbrace\begin{bmatrix}
a & b\\ c& d
\end{bmatrix} \in \mbox{SL}_2(\mathbb{Z}) : a-1\equiv c\equiv d-1\equiv 0\pmod N\right\rbrace,\\
\Gamma_0(N) &:= \left\lbrace\begin{bmatrix}
a & b\\ c& d
\end{bmatrix} \in \mbox{SL}_2(\mathbb{Z}) : c\equiv 0\pmod N\right\rbrace.
\end{aligned}\]
Any element of a given congruence subgroup $\Gamma$ acts on the extended upper half-plane $\mathbb{H}^\ast := \mathbb{H}\cup \mathbb{Q}\cup \{\infty\}$ by a linear fractional transformation. We define the cusps of $\Gamma$ as the equivalence classes of $\mathbb{Q}\cup \{\infty\}$ under this action, and it is known that there are finitely many inequivalent cusps of $\Gamma$ (see \cite[Prop. 1.32]{shim}). We define the modular function on $\Gamma$ by a function $f(\tau):\mathbb{H}\rightarrow\mathbb{C}$ such that
\begin{itemize}
\item $f(\tau)$ is meromorphic on $\mathbb{H}$,
\item $f(\gamma\tau)=f(\tau)$ for all $\gamma\in \Gamma$, and 
\item $f(\tau)$ is meromorphic at all cusps of $\Gamma$.
\end{itemize}
The last condition implies that for every cusp $r$ of $\Gamma$ and an element $\gamma\in\mbox{SL}_2(\mathbb{Z})$ with $\gamma(\infty)=r$, the $q$-expansion of $f\circ \gamma$ given by $f(\gamma\tau) = \sum_{n\geq n_0} a_nq^{n/h}$, for some integers $h$ and $n_0$ with $a_{n_0}\neq 0$, has finitely many negative exponents. We call the least such $n_0$ the order of $f(\tau)$ at $r$, denoted by $\mbox{ord}_r f(\tau)$, and we say that $f(\tau)$ has a zero (resp., a pole) at $r$ if $\mbox{ord}_r f(\tau)$ is positive (resp., negative). We also call $h$ the width of $r$, which is the smallest positive integer such that $\gamma[\begin{smallmatrix}
1 & h\\0 & 1
\end{smallmatrix}]\gamma^{-1}\in \pm\Gamma$. We note that $h$ is independent on the choice of $\gamma$.

Let $A_0(\Gamma)$ be the field of all modular functions on $\Gamma$ and let $A_0(\Gamma)_{\mathbb{Q}}$ be its subfield consisting of all modular functions on $\Gamma$ with rational $q$-expansion. We identify $A_0(\Gamma)$ with the field of all meromorphic functions on the modular curve $X(\Gamma):=\Gamma\backslash \mathbb{H}^\ast$. If $f(\tau)\in A_0(\Gamma)$ has no zeros nor poles on $\mathbb{H}$, then the degree $[A_0(\Gamma):\mathbb{C}(f(\tau))]$ is the total degree of poles of $f(\tau)$ given by $-\sum_r \mbox{ord}_r f(\tau)$, where the sum ranges over the inequivalent cusps $r$ of $\Gamma$ for which $f(\tau)$ has a pole at $r$ (see \cite[Prop. 2.11]{shim}).

We describe a systematic way of getting the set of inequivalent cusps of the congruence subgroup $\Gamma:=\Gamma_1(N)\cap \Gamma_0(mN)$ for positive integers $m$ and $N$, following Cho, Koo, and Park \cite{chokoopark}. Let $\Gamma(1)_{\infty}$ be the stabilizer of $\infty$ in $\mbox{SL}_2(\mathbb{Z})$ and let $\Gamma\gamma_i\Gamma(1)_{\infty}, i=1,\ldots,g,$ be the representatives
for the double coset $\Gamma\backslash \mbox{SL}_2(\mathbb{Z})/\Gamma(1)_{\infty}$. Then the set $\{\gamma_i(\infty):i=1,\ldots,g\}$ is the set of all inequivalent cusps of $\Gamma$. Consider the set
\begin{align*}
M := \{ (c,d)\in (\mathbb{Z}/mN\mathbb{Z})^2 : \gcd(c,d,mN) = 1\}
\end{align*}
and the subgroup $\Delta$ of $(\mathbb{Z}/mN\mathbb{Z})^\times$ given by
\begin{align*}
\Delta := \{\pm(1+Nk)\in (\mathbb{Z}/mN\mathbb{Z})^\times : k=0,\ldots,m-1\}.
\end{align*}
We define a relation $\sim$ on $M$ as follows: if $(c_1,d_1), (c_2,d_2)\in M$, then $(c_1,d_1)\sim (c_2,d_2)$ if $(c_2,d_2)\equiv(sc_1,sd_1+nc_1)\pmod{mN}$ for some $s\in \Delta$ and $n\in \mathbb{Z}/mN\mathbb{Z}$. This is an equivalence relation, and the map $\phi: \Gamma\backslash\mbox{SL}_2(\mathbb{Z})/\Gamma(1)_{\infty}\rightarrow M/\sim$ sending $\Gamma[\begin{smallmatrix}
a&b\\ c& d\end{smallmatrix}]\Gamma(1)_{\infty}$ to the equivalence class of $(c,d)$ is well-defined and bijective (see \cite[Prop. 3.8.5]{diashur}).

\begin{lemma}[\cite{chokoopark}]\label{lem21}
Suppose $a,c,a',c'\in\mathbb{Z}$ such that $\gcd(a,c)=\gcd(a',c')=1$. We denote $\infty := \pm 1/0$. Then $a/c$ and $a'/c'$ are equivalent under $\Gamma_1(N)\cap \Gamma_0(mN)$ if and only if there exist $s\in\Delta$ and $n\in\mathbb{Z}$ such that $[\begin{smallmatrix}
a' \\ c'
\end{smallmatrix}]\equiv [\begin{smallmatrix}
s^{-1}a+nc \\ sc
\end{smallmatrix}]\pmod{mN}$.
\end{lemma}

We now consider the surjective canonical projection $\pi_x : (\mathbb{Z}/mN\mathbb{Z})^\times\rightarrow (\mathbb{Z}/x\mathbb{Z})^\times$ for a positive divisor $x$ of $mN$. For each positive divisor $c$ of $mN$, let $s'_{c,i}, i=1,\ldots, n_c$, be the distinct coset representatives of $\pi_{mN/c}(\Delta)$ in $(\mathbb{Z}/\frac{mN}{c}\mathbb{Z})^\times$ with $n_c = \varphi(\frac{mN}{c})/|\pi_{mN/c}(\Delta)|$ (here, $\varphi$ is the Euler's totient function). We then choose $s_{c,i}\in(\mathbb{Z}/mN\mathbb{Z})^\times$ such that $\pi_{mN/c}(s_{c,i})=s'_{c,i}$ for each $i$ and define $S_c := \{s_{c,i}\in(\mathbb{Z}/mN\mathbb{Z})^\times: i=1,\ldots,n_c\}$. We next let $a'_{c,j}, j=1,\ldots, m_c$, be the distinct coset representatives of $\pi_c(\Delta\cap \ker \pi_{mN/c})$ in $(\mathbb{Z}/c\mathbb{Z})^\times$ with 
\begin{align*}
m_c = \dfrac{\varphi(c)}{|\pi_c(\Delta\cap \ker \pi_{mN/c})|} = \dfrac{\varphi(c)|\pi_{mN/c}(\Delta)|}{|\pi_{mN/\gcd(c,mN/c)}(\Delta)|}.
\end{align*}
We then choose a positive integer $a_{c,j}$ less than $mN$ and coprime to $mN$ such that $\pi_{c}(a_{c,j})=a'_{c,j}$ for each $j$ and define $A_c := \{a_{c,j}: j=1,\ldots,m_c\}$. 

\begin{lemma}[\cite{chokoopark}]\label{lem22}
With the notations as above, let
{\small \begin{align*}
S=\{(cs_{c,i},a_{c,j})\in  (\mathbb{Z}/mN\mathbb{Z})^2 : 0 < c\mid mN, s_{c,i}\in S_c, a_{c,j}\in A_c, \gcd(cs_{c,i},a_{c,j},mN)=1\}.
\end{align*} }%
For a given $(cs_{c,i},a_{c,j})\in S$, choose coprime integers $x$ and $y$ such that $x\equiv cs_{c,i}\pmod{mN}$ and $y\equiv a_{c,j}\pmod{mN}$. Then the set of $y/x$ for such $x$ and $y$ is the set of all inequivalent cusps of $\Gamma_1(N)\cap \Gamma_0(mN)$ and the number of such cusps is 
\begin{align*}
|S| = \sum_{0<c\mid mN} n_cm_c = \sum_{0<c\mid mN} \dfrac{\varphi(c)\varphi(mN/c)}{|\pi_{mN/\gcd(c,mN/c)}(\Delta)|}.
\end{align*}
\end{lemma}

We now obtain the width of a cusp of $\Gamma_1(N)\cap \Gamma_0(mN)$ as the following result shows.

\begin{lemma}[\cite{chokoopark}]\label{lem23}
The width $h$ of a cusp $a/c$ of $\Gamma_1(N)\cap \Gamma_0(mN)$, where $a$ and $c$ are coprime integers, is  
\begin{align*}
h &= \begin{cases}
\dfrac{m}{\gcd(c^2/4,m)} &\text{ if }N=4, \gcd(m,2)=1,\gcd(c,4)=2,\\
\dfrac{mN}{\gcd(c,N)\gcd(m,c^2/\gcd(c,N))} &\text{ otherwise. }
\end{cases}
\end{align*}
\end{lemma}

For the congruence subgroups $\Gamma_0(N)$ and $\Gamma_1(N)$, we obtain a nice description for their respective sets of inequivalent cusps.

\begin{corollary}[\cite{chokoopark}]\label{cor24}
Suppose $a,c,a',c'\in\mathbb{Z}$ such that $\gcd(a,c)=\gcd(a',c')=1$. We denote $\infty := \pm 1/0$. Let $S_\Gamma$ be the set of all inequivalent cusps of a congruence subgroup $\Gamma$.
\begin{enumerate}
\item[(1)] The cusps $a/c$ and $a'/c'$ are equivalent under $\Gamma_0(m)$ if and only if there exist $s\in(\mathbb{Z}/m\mathbb{Z})^\times$ and $n\in\mathbb{Z}$ such that $[\begin{smallmatrix}
a' \\ c'
\end{smallmatrix}]\equiv [\begin{smallmatrix}
s^{-1}a+nc \\ sc
\end{smallmatrix}]\pmod{m}$. Moreover, we have 
\begin{align*}
S_{\Gamma_0(m)} =&\{a_{c,j}/c\in\mathbb{Q} : 0 < c\mid m, 0<a_{c,j}\leq m, \gcd(a_{c,j},m)=1\\
& a_{c,j}=a_{c,j'}\iff a_{c,j}\equiv a_{c,j'}\pmod{\gcd(c,m/c)}\}
\end{align*}
and the width of the cusp $a/c$ in $\Gamma_0(m)$ is $m/\gcd(c^2,m)$.
\item[(2)] The cusps $a/c$ and $a'/c'$ are equivalent under $\Gamma_1(N)$ if and only if there exists $n\in\mathbb{Z}$ such that $[\begin{smallmatrix}
a' \\ c'
\end{smallmatrix}]\equiv \pm[\begin{smallmatrix}
a+nc \\ c
\end{smallmatrix}]\pmod{N}$. Moreover, we have 
\begin{align*}
S_{\Gamma_1(N)} =&\{y_{c,j}/x_{c,i}\in\mathbb{Q} : 0 < c\mid N, 0<s_{c,i},a_{c,j}\leq N,\\ 
&\gcd(s_{c,i},N)=\gcd(a_{c,j},N)=1,\\
& s_{c,i}=s_{c,i'}\iff s_{c,i}\equiv \pm s_{c,i'}\pmod{N/c},\\
& a_{c,j}=a_{c,j'}\iff a_{c,j}\equiv \begin{cases}
	\pm a_{c,j'}\pmod{c} &\text{ if } c\in\{N/2,N\},\\
	a_{c,j'}\pmod{c} &\text{ otherwise, }
\end{cases}\\
&(x_{c,i},y_{c,j})\in\mathbb{Z}^2, \gcd(x_{c,i},y_{c,j})=1,\\
&x_{c,i}\equiv cs_{c,i}\pmod{N}, y_{c,j}\equiv a_{c,j}\pmod{N}\}
\end{align*}
and the width of the cusp $a/c$ in $\Gamma_1(N)$ is $1$ if $N=4$ and $\gcd(c,4)=2$, and $N/\gcd(c,N)$ otherwise.
\end{enumerate}
\end{corollary}	

We now introduce an $\eta$-quotient, which is a function of the form
\begin{equation*}
f(\tau) = \prod_{\delta\mid N} \eta(\delta\tau)^{r_{\delta}}
\end{equation*}
for some indexed set $\{r_\delta\in\mathbb{Z} : \delta\mid N\}$, where $\eta(\tau)=q^{1/24}\prod_{n=1}^\infty(1-q^n)$ is the Dedekind eta function. The next two lemmas provide necessary conditions on the modularity and the behavior at the cusps of an $\eta$-quotient on $\Gamma_0(N)$.

\begin{lemma}[\cite{newman1, newman2, ono}]\label{lem25}
Let $f(\tau) = \prod_{\delta\mid N} \eta(\delta\tau)^{r_{\delta}}$ be an $\eta$-quotient with $k= \frac{1}{2}\sum_{\delta\mid N} r_{\delta}\in\mathbb{Z}$ such that
\begin{equation*}
\sum_{\delta\mid N} \delta r_{\delta} \equiv 0\pmod{24}\quad\text{ and }\quad\sum_{\delta\mid N} \dfrac{N}{\delta}r_{\delta} \equiv 0\pmod {24}.
\end{equation*}
Then for all $[\begin{smallmatrix}
a & b\\ c & d
\end{smallmatrix}]\in \Gamma_0(N)$,
\begin{equation*}
f\left(\dfrac{a\tau+b}{c\tau+d}\right) = \left(\dfrac{(-1)^k\prod_{\delta\mid N} \delta^{r_\delta}}{d}\right)(c\tau+d)^kf(\tau),
\end{equation*}
where $(\frac{\cdot}{d})$ is the Kronecker symbol.
\end{lemma} 

\begin{lemma}[\cite{ligo,ono}]\label{lem26}
Let $c, d$ and $N$ be positive integers with $d\mid N$ and $\gcd(c,d)=1$ and let $f(\tau) = \prod_{\delta\mid N} \eta(\delta\tau)^{r_{\delta}}$ be an $\eta$-quotient satisfying the conditions of Lemma \ref{lem25}. Then the order of vanishing of $f(\tau)$ at the cusp $c/d$ is 
\begin{align*}
\dfrac{N}{24d\gcd(d,\frac{N}{d})}\sum_{\delta\mid N} \gcd(d,\delta)^2\cdot\dfrac{r_{\delta}}{\delta}.
\end{align*}
\end{lemma}

We next define, following Yang \cite{yang}, generalized $\eta$-quotients with analogous properties on $\Gamma_1(N)$. These are functions of the form
\begin{equation*}
f(\tau) = \prod_{1\leq g\leq \lfloor N/2\rfloor}\eta_{N,g}(\tau)^{r_g}
\end{equation*}
for some indexed set $\{r_g \in\mathbb{Z}: 1\leq g\leq \lfloor N/2\rfloor\}$ of integers, where 
\begin{equation*}
\eta_{N, g}(\tau) = q^{NB_2(g/N)/2}\prod_{m=1}^\infty(1-q^{N(m-1)+g})(1-q^{Nm-g})
\end{equation*}
is the generalized Dedekind eta function with $B_2(t) := t^2-t+\frac{1}{6}$. The following results give the transformation formula of $\eta_{N, g}$ on $\Gamma_0(N)$ and conditions for the modularity of a generalized $\eta$-quotient on $\Gamma_1(N)$.

\begin{lemma}[\cite{yang}]\label{lem27}
The function $\eta_{N,g}(\tau)$ satisfies $\eta_{N,g+N}(\tau)=\eta_{N,-g}(\tau)=-\eta_{N,g}(\tau)$. Moreover, let $\gamma = [\begin{smallmatrix}
a & b\\ cN & d
\end{smallmatrix}]\in \Gamma_0(N)$ with $c\neq 0$. Then we have 
\begin{align*}
\eta_{N,g}(\gamma\tau) = \varepsilon(a,bN,c,d)e^{\pi i(g^2ab/N - gb)}\eta_{N,ag}(\tau),
\end{align*}
where 
\begin{align*}
\varepsilon(a,b,c,d) = \begin{cases}
e^{\pi i(bd(1-c^2)+c(a+d-3))/6}, & \text{ if }c\equiv 1\pmod{2},\\
-ie^{\pi i(ac(1-d^2)+d(b-c+3))/6}, & \text{ if }c\equiv 0\pmod{2}.
\end{cases}
\end{align*}
\end{lemma}

\begin{lemma}[\cite{yang}]\label{lem28}
Suppose $f(\tau)=\prod_{1\leq g\leq \lfloor N/2\rfloor}\eta_{N,g}(\tau)^{r_g}$ is a generalized $\eta$-quotient such that
\begin{align*}
\sum_{1\leq g\leq \lfloor N/2\rfloor}r_g\equiv 0\pmod{12}\quad\text{ and }\quad\sum_{1\leq g\leq \lfloor N/2\rfloor}gr_g\equiv 0\pmod{2}.
\end{align*}
Then $f(\tau)$ is a modular function on $\Gamma(N)$. Moreover, if in addition, $f(\tau)$ satisfies
\begin{align*}
\sum_{1\leq g\leq \lfloor N/2\rfloor}g^2r_g\equiv 0\pmod{2N}
\end{align*}
then $f(\tau)$ is a modular function on $\Gamma_1(N)$. 
\end{lemma}

\begin{lemma}[\cite{yang}]\label{lem29}
Let $N$ be a positive integer and $\gamma=[\begin{smallmatrix}
a & b\\ c& d
\end{smallmatrix}]\in\mathrm{SL}_2(\mathbb{Z})$. Then the first term of the $q$-expansion of $\eta_{N, g}(\gamma\tau)$ is $\varepsilon q^{\delta}$, where $|\varepsilon|=1$ and 
\[\delta = \dfrac{\gcd(c,N)^2}{2N}P_2\left(\dfrac{ag}{\gcd(c,N)}\right)\]
with $P_2(t):=B_2(\{t\})$ and $\{t\}$ is the fractional part of $t$.
\end{lemma}

Using Lemma \ref{lem23}, Corollary \ref{cor24}, and Lemma \ref{lem29}, we can compute the order of a generalized $\eta$-quotient at a cusp of $\Gamma_1(N)\cap \Gamma_0(mN)$.	

\section{Modularity of \texorpdfstring{$r(\tau)^ar(2\tau)^b$}{r(t)ar(t)b}}\label{sec3}

In this section, we show that under a congruence condition on integers $a$ and $b$ with $(a,b)\neq (0,0)$, there are finitely many functions of the form $r(\tau)^ar(2\tau)^b$ generating the field $A_0(\Gamma_1(10))$ of all modular functions on $\Gamma_1(10)$, partially extending the result of Lee and Park \cite{leepark} that $k(\tau)$ is one such function. To prove this, we use the fact that $r(\tau) = \eta_{5, 1}(\tau)/\eta_{5, 2}(\tau)$.

\begin{theorem}\label{thm31}
Let $a$ and $b$ be integers with $(a,b)\neq (0,0)$ such that $a\equiv 3b\pmod{5}$. Then the function $r(\tau)^ar(2\tau)^b$ is modular on $\Gamma_1(10)$. In addition, if $r(\tau)^ar(2\tau)^b$ generates $A_0(\Gamma_1(10))$, then there are two such functions up to taking reciprocals, namely $k(\tau)$ and $l(\tau)$.
\end{theorem}

\begin{proof}
We write $F(\tau):=r(\tau)^ar(2\tau)^b$ as a generalized $\eta$-quotient given by 
\begin{align*}
F(\tau) = \eta_{10, 1}(\tau)^a\eta_{10, 2}(\tau)^{b-a}\eta_{10, 3}(\tau)^{-a}\eta_{10, 4}(\tau)^{a-b}.
\end{align*}
By Lemma \ref{lem28} and the given congruence condition, we see that $F(\tau)$ is modular on $\Gamma_1(10)$. Using Corollary \ref{cor24}(2), we choose the set of inequivalent cusps of $\Gamma_1(10)$ to be $S_{\Gamma_1(10)}=\{\infty, 0,1/2,1/3,1/5,1/6,3/5,3/10\}$. We then calculate using Lemma \ref{lem29} the order of $F(\tau)$ at each element of $S_{\Gamma_1(10)}$, as shown in Table \ref{tab:tbl1}.

\begin{table}[h]
\begin{tabular}{@{}cccccc@{}}
\toprule
cusp $r$ & $\infty$  & $1/5$  & $3/5$ & $3/10$ & $0,1/2,1/3,1/6$ \\
\midrule
$\mbox{ord}_r F(\tau)$ & $\frac{1}{5}(a+2b)$ & $\frac{1}{5}(2a-b)$ & $\frac{1}{5}(-2a+b)$ &  $\frac{1}{5}(-a-2b)$ & $0$ \\
\bottomrule
\end{tabular}
\caption{The orders of $F(\tau)$ at the cusps of $\Gamma_1(10)$}\label{tab:tbl1}
\end{table}

Observe that the orders of $F(\tau)$ at the cusps $\infty, 1/5, 3/5$ and $3/10$ are integers, and that $\mbox{ord}_{3/5} F(\tau)=- \mbox{ord}_{1/5} F(\tau)$ and $\mbox{ord}_{3/10} F(\tau)=- \mbox{ord}_\infty F(\tau)$. Since $F(\tau)$ generates $A_0(\Gamma_1(10))$, exactly one of $\mbox{ord}_{1/5} F(\tau)$ and $\mbox{ord}_\infty F(\tau)$ is zero and the other is $\pm 1$ (otherwise, the total degree of poles of $F(\tau)$ would have been at least $2$). We now proceed by doing casework:	
\begin{itemize}
\item If $0=\mbox{ord}_{1/5} F(\tau)=\frac{1}{5}(2a-b)$ and $\pm 1=\mbox{ord}_\infty F(\tau)=\frac{1}{5}(a+2b)$, then $(a,b)\in\{(1,2),(-1,-2)\}$.
\item If $0=\mbox{ord}_\infty F(\tau)=\frac{1}{5}(a+2b)$ and $\pm 1=\mbox{ord}_{1/5} F(\tau)=\frac{1}{5}(2a-b)$, then $(a,b)\in\{(-2,1),(2,-1)\}$.
\end{itemize}
Thus, we get four possible ordered pairs $(a,b)$, corresponding to the functions $k(\tau), k(\tau)^{-1}, l(\tau)$ and $l(\tau)^{-1}$. The desired conclusion now follows.
\end{proof}

As a consequence of Theorem \ref{thm31}, we prove that the field $A_0(\Gamma_0(10))$ of all modular functions on $\Gamma_0(10)$ can be generated by a rational function of $l(\tau)$. We use the fact that $\Gamma_0(10)$ is generated by $\Gamma_1(10)$ and the matrix $\gamma:=[\begin{smallmatrix}
3 & -1\\10 & -3
\end{smallmatrix}]$.

\begin{theorem}\label{thm32}
The function $l(\tau)-1/l(\tau)$ generates $A_0(\Gamma_0(10))$.
\end{theorem}

\begin{proof}
Since $l(\gamma\tau) = -1/l(\tau)$ by Lemma \ref{lem27}, $l(\tau)-1/l(\tau)\in A_0(\Gamma_0(10))$. We take $S_{\Gamma_0(10)} = \{\infty,0,1/2,1/5\}$ as the 
set of inequivalent cusps of $\Gamma_0(10)$ by Corollary \ref{cor24}(1). Using Lemma \ref{lem29}, we compute $\mbox{ord}_{1/5} l(\tau)=-1$ and $\mbox{ord}_r l(\tau)=0$ for any $r\in S_{\Gamma_0(10)}$ with $r\neq 1/5$. We see that 
\begin{align*}
\mbox{ord}_{1/5} (l(\tau)-1/l(\tau)) = \min\{\mbox{ord}_{1/5} l(\tau), \mbox{ord}_{1/5} 1/l(\tau)\}=-1
\end{align*}
and 
\begin{align*}
\mbox{ord}_r (l(\tau)-1/l(\tau)) \geq \min\{\mbox{ord}_r l(\tau), \mbox{ord}_r 1/l(\tau)\}=0
\end{align*}
for any $r\in S_{\Gamma_0(10)}$ with $r\neq 1/5$. Thus, we obtain that $l(\tau)-1/l(\tau)$ has a simple pole only at $1/5$ and holomorphic elsewhere, which yields the desired conclusion.
\end{proof}

We now consider modular functions on $\Gamma_0(10)$ involving $r(\tau)^ar(2\tau)^b$ for integers $a$ and $b$ with $(a,b)\neq (0,0)$. The following results, which extend Theorem \ref{thm32}, imply that
under the condition of Theorem \ref{thm31}, the functions $r(\tau)^ar(2\tau)^b$ can be expressed explicitly in terms of $l(\tau)$.

\begin{theorem}\label{thm33}
Let $a$ and $b$ be integers with $(a,b)\neq (0,0)$ such that $a\equiv 3b\pmod{5}$. Then the function $f_{a,b}(\tau) := r(\tau)^ar(2\tau)^b+(-1)^{b-a}/(r(\tau)^ar(2\tau)^b)$ is modular on $\Gamma_0(10)$. In addition, the function $(l(\tau)-1/l(\tau))^{|a+2b|/5}f_{a,b}(\tau)$ is a polynomial in $l(\tau)-1/l(\tau)$ with rational coefficients.
\end{theorem}

\begin{proof}
We know from Theorem \ref{thm31} that $F(\tau):=r(\tau)^ar(2\tau)^b$ is a modular function on $\Gamma_1(10)$. Since $F(\gamma\tau) = (-1)^{b-a}/F(\tau)$ by Lemma \ref{lem27}, we see that $f_{a,b}(\tau) = F(\tau)+F(\gamma\tau)$ is a modular function on $\Gamma_0(10)$. We again choose $S_{\Gamma_0(10)} = \{\infty,0,1/2,1/5\}$ as the 
set of inequivalent cusps of $\Gamma_0(10)$ by Corollary \ref{cor24}(1). By Lemma \ref{lem29}, we find that $\mbox{ord}_{\infty} F(\tau)=(a+2b)/5, \mbox{ord}_{1/5} F(\tau)=(2a-b)/5$ and $\mbox{ord}_r F(\tau)=0$ for $r\in \{0,1/2\}$. Thus, we get
\begin{align*}
\mbox{ord}_{\infty} f_{a,b}(\tau) &=  \min\{\mbox{ord}_{\infty} F(\tau), \mbox{ord}_{\infty} 1/F(\tau)\}=-|a+2b|/5,\\
\mbox{ord}_{1/5} f_{a,b}(\tau) &=  \min\{\mbox{ord}_{1/5} F(\tau), \mbox{ord}_{1/5} 1/F(\tau)\}=-|2a-b|/5,
\end{align*}
and 
\begin{align*}
\mbox{ord}_r f_{a,b}(\tau) \geq  \min\{\mbox{ord}_r F(\tau), \mbox{ord}_r 1/F(\tau)\}=0
\end{align*}
for $r\in \{0,1/2\}$. We now claim that $l(\tau)-1/l(\tau)$ has a zero at $\infty$. Indeed, we have the $q$-expansion of $l(\tau)$ given by
\begin{equation}\label{eq31}
l(\tau) = 1 + 2 q - 4 q^3 - 2 q^4 + 6 q^5 + 8 q^6 - 4 q^7 +O(q^8).
\end{equation}
Observe that when $\tau=\infty$, $q=0$ and $l(\tau)=1$, which proves our claim. Consequently, we see that $\mbox{ord}_{\infty}(l(\tau)-1/l(\tau))\geq 1$. We deduce from Theorem \ref{thm32} that the
nonconstant function $(l(\tau)-1/l(\tau))^{|a+2b|/5}f_{a,b}(\tau)$ has a pole of order $d:=(|a+2b|+|2a-b|)/5$ only at $1/5$ and holomorphic elsewhere. Thus, by \cite[Lem. 2]{yang2}, we obtain that $(l(\tau)-1/l(\tau))^{|a+2b|/5}f_{a,b}(\tau)$ is a polynomial in $l(\tau)-1/l(\tau)$ of degree $d$. Moreover, this polynomial has rational coefficients since the $q$-expansions of $r(\tau)$ and $l(\tau)$ are rational.
\end{proof}

\begin{corollary}\label{cor34}
We have the identities
\begin{align*}
r(\tau)^5 = l^{-2}\left(\dfrac{l-1}{l+1}\right), \qquad r(2\tau)^5 = l\left(\dfrac{l-1}{l+1}\right)^2,
\end{align*}
where $l:=l(\tau)$. Consequently, under the conditions of Theorem \ref{thm33}, we have
\begin{align*}
r(\tau)^ar(2\tau)^b = l^{(b-2a)/5}\left(\dfrac{l-1}{l+1}\right)^{(a+2b)/5}.
\end{align*}
\end{corollary}

\begin{proof}
Using the notations of Theorem \ref{thm33} and following its proof, we know that $(l-1/l)f_{5,0}(\tau)$ is a third-degree polynomial in $l-1/l$. Writing
\begin{align*}
(l-1/l)f_{5,0}(\tau) = \sum_{i=0}^3 a_i(l-1/l)^i
\end{align*}
for some constants $a_0,a_1,a_2$ and $a_3$, and applying the $q$-expansions
\begin{align*}
r(\tau)^5 =q - 5 q^2 + 15 q^3 - 30 q^4 + 40 q^5 - 26 q^6+O(q^7)
\end{align*}
and (\ref{eq31}), we get $(a_0,a_1,a_2,a_3)=(-4,-4,-2,-1)$ and 
\begin{align*}
f_{5,0}(\tau)=r(\tau)^5-\dfrac{1}{r(\tau)^5}=-\dfrac{(1-l+l^2+l^3)(-1+l+l^2+l^3)}{l^2(-1+l)(1+l)}.
\end{align*}
Solving the above equation for $r(\tau)^5$ and applying again the $q$-expansions of $r(\tau)^5$ and $l(\tau)$ yield the first identity. The identities involving $r(2\tau)^5$ and $r(\tau)^ar(2\tau)^b$ now follow from the first identity and the definition of $l(\tau)$.
\end{proof}

\section{Modular equations for \texorpdfstring{$l(\tau)$}{l(t)}}\label{sec4}

We infer from Theorem \ref{thm31} that $l(\tau)$ generates $A_0(\Gamma_1(10))$. The following shows that we can define an affine model for the modular curve $X(\Gamma_1(10)\cap \Gamma_0(10n))$ over $\mathbb{Q}$ using $l(\tau)$. This affine model gives rise to the modular equation for $l(\tau)$ of suitable level. We omit the proof as it is similar to that of \cite[Lem. 3.2]{leepark}. 

\begin{lemma}\label{lem41}
We have $A_0(\Gamma_1(10)\cap \Gamma_0(10n))_{\mathbb{Q}}=\mathbb{Q}(l(\tau),l(n\tau))$ for all positive integers $n$.
\end{lemma}

According to Table \ref{tab:tbl1}, $l(\tau)$ has a simple pole at $1/5$ and a simple zero at $3/5$. We use this to describe the behavior of $l(\tau)$ at the cusps.  

\begin{lemma}\label{lem42}
Let $a, c, a', c'$ and $n$ be integers with $n$ positive. Then:
\begin{enumerate}
\item[(1)] The modular function $l(\tau)$ has a pole at $a/c\in \mathbb{Q}\cup\{\infty\}$ if and only if $\gcd(a,c) = 1, a\equiv \pm 1\pmod{5}$ and $c\equiv 5\pmod {10}$.
\item[(2)] The modular function $l(n\tau)$ has a pole at $a'/c'\in \mathbb{Q}\cup\{\infty\}$ if and only if there exist integers $a$ and $c$ such that $a/c=na'/c', \gcd(a,c)=1, a\equiv \pm 1\pmod{5}$ and $c\equiv 5\pmod{10}$.
\item[(3)] The modular function $l(\tau)$ has a zero at $a/c\in \mathbb{Q}\cup\{\infty\}$ if and only if $\gcd(a,c) = 1, a\equiv \pm 3\pmod{5}$ and $c\equiv 5\pmod {10}$.
\item[(4)] The modular function $l(n\tau)$ has a zero at $a'/c'\in \mathbb{Q}\cup\{\infty\}$ if and only if there exist integers $a$ and $c$ such that $a/c=na'/c', \gcd(a,c)=1, a\equiv \pm 3\pmod{5}$ and $c\equiv 5\pmod{10}$.
\end{enumerate}
\end{lemma}

\begin{proof}
Note that $l(\tau)$ has a pole (resp., zero) at $a/c\in \mathbb{Q}\cup\{\infty\}$ which is equivalent to $1/5$ (resp., $3/5$) on $\Gamma_1(10)$. By Corollary \ref{cor24}(2), $a/c$ is equivalent to $1/5$ on $\Gamma_1(10)$ if and only if $[\begin{smallmatrix}
a \\ c
\end{smallmatrix}]\equiv \pm[\begin{smallmatrix}
1+5m \\ 5
\end{smallmatrix}]\pmod{10}$ for some integer $m$. Thus, we get $\gcd(a,c)=1, a\equiv \pm 1\pmod{5}$ and $c\equiv 5\pmod{10}$, proving (1). On the other hand, $a/c$ is equivalent to $3/5$ on $\Gamma_1(10)$ if and only if $[\begin{smallmatrix}
a \\ c
\end{smallmatrix}]\equiv \pm[\begin{smallmatrix}
3+5m \\ 5
\end{smallmatrix}]\pmod{10}$ for some integer $m$. This yields $\gcd(a,c)=1, a\equiv \pm 3\pmod{5}$ and $c\equiv 5\pmod{10}$, proving (3). Statement (2) (resp. (4)) now follows from (1)  (resp. (3)).
\end{proof}

To obtain explicitly modular equations for $l(\tau)$ of any level, we present the following result by Ishida and Ishii \cite{isdishii}, which gives information on the coefficients of the modular equation between two functions $f_1(\tau)$ and $f_2(\tau)$ that generate the field $A_0(\Gamma')$ for some congruence subgroup $\Gamma'$ of $\mbox{SL}_2(\mathbb{Z})$.

\begin{proposition}[\cite{chokimkoo,isdishii}]\label{prop43}
Let $\Gamma'$ be a congruence subgroup of $\mathrm{SL}_2(\mathbb{Z})$, and $f_1(\tau)$ and $f_2(\tau)$ be nonconstant functions such that $A_0(\Gamma')=\mathbb{C}(f_1(\tau),f_2(\tau))$. For $k\in\{1,2\}$, let $d_k$ be the total degree of poles of $f_k(\tau)$. Let 
\begin{align*}
F(X,Y)=\sum_{\substack{0\leq i\leq d_2\\0\leq j\leq d_1}}C_{i,j}X^iY^j\in \mathbb{C}[X,Y]
\end{align*}
satisfy $F(f_1(\tau),f_2(\tau))=0$. Let $S_{\Gamma'}$ be the set of all the inequivalent cusps of $\Gamma'$, and for $k\in\{1,2\}$, define 
\begin{align*}
S_{k,0}&= \{r\in S_{\Gamma'} : f_k(\tau)\text{ has a zero at }r\},\\
S_{k,\infty}&= \{r\in S_{\Gamma'} : f_k(\tau)\text{ has a pole at }r\}.
\end{align*}
Further let
\begin{align*}
a = -\sum_{r\in S_{1,\infty}\cap S_{2,0}}\mbox{ord}_r f_1(\tau)\quad\text{ and }\quad b = \sum_{r\in S_{1,0}\cap S_{2,0}}\mbox{ord}_r f_1(\tau)
\end{align*}
with the convention that $a=0$ (resp. $b=0$) whenever $S_{1,\infty}\cap S_{2,0}$ (resp. $S_{1,0}\cap S_{2,0}$) is empty. Then
\begin{enumerate}
\item[(1)] $C_{d_2,a}\neq 0$ and if, in addition, $S_{1,\infty}\subseteq S_{2,\infty}\cup S_{2,0}$, then $C_{d_2,j}=0$ for all $j\neq a$;
\item[(2)] $C_{0,b}\neq 0$ and if, in addition, $S_{1,0}\subseteq S_{2,\infty}\cup S_{2,0}$, then $C_{0,j}=0$ for all $j\neq b$.
\end{enumerate}
By interchanging the roles of $f_1(\tau)$ and $f_2(\tau)$, one may obtain properties analogous to (1) - (2).
\end{proposition}

We apply Lemma \ref{lem42} and Proposition \ref{prop43} to deduce the following result, which asserts that modular equations for $l(\tau)$ of any level always exist, which serve as affine models for $X(\Gamma_1(10)\cap \Gamma_0(10n))$.

\begin{theorem}\label{thm44}
One can explicitly obtain modular equation for $l(\tau)$ of level $n$ for all positive integers $n$.
\end{theorem}

\begin{proof}
The proof is similar to that of \cite[Thm. 1.1]{leepark}(1), so we omit the details. We instead note that if $L_n(l(\tau),l(n\tau))=0$ is the modular equation for $l(\tau)$ of level $n$ for some polynomial $L_n(X,Y)\in\mathbb{C}[X,Y]$, then $L_n(X,Y)$ is irreducible as a polynomial of $X$ (resp. $Y$) over $\mathbb{C}(Y)$ (resp. $\mathbb{C}(X)$), in view of \cite[Lem. 3.1]{chokimkoo}.
\end{proof}

We now use Theorem \ref{thm44} to derive the modular equations for $l(\tau)$ of levels two and five, and find which coefficients of the modular equations for $l(\tau)$ of odd prime level $p\neq 5$ are nonzero.

\begin{theorem}[Modular equation of level two]\label{thm45}
We have
\begin{align*}
1+l(\tau)-2l(\tau)l(2\tau)-l(\tau)l(2\tau)^2+l(\tau)^2l(2\tau)^2=0.
\end{align*}
\end{theorem}

\begin{proof}
In view of Lemmas \ref{lem22} and \ref{lem42}, we work on the congruence subgroup $\Gamma_1(10)\cap \Gamma_0(20)$ and the following cusps: $1/5, 3/5, 1/10$ and $3/10$. Recall that 
\[l(\tau) = \dfrac{\eta_{10, 2}(\tau)^3\eta_{10, 3}(\tau)^2}{\eta_{10, 1}(\tau)^2\eta_{10, 4}(\tau)^3}\]
is a generalized $\eta$-quotient. Using Lemmas \ref{lem23} and \ref{lem29}, we know that $l(\tau)$ has a double pole at $1/5$ and a double zero at $3/5$. We now compute the order of $l(2\tau)$ at $3/10$ (which is a zero) as follows. Since 
\[l\left(2\begin{bmatrix}
3 & -1\\10 & -3
\end{bmatrix}\tau\right)=l\left(\begin{bmatrix}
6 & -2\\10 & -3
\end{bmatrix}\tau\right)=l\left(\begin{bmatrix}
3 & -2\\5 & -3
\end{bmatrix}(2\tau)\right),\]
by Lemma \ref{lem29}, the first term of the $q$-expansion of $l(2[\begin{smallmatrix}
3 & -1\\10 & -3
\end{smallmatrix}]\tau)$ is $\varepsilon q^\delta$, where $|\varepsilon|=1$ and 
\[\delta = 2\cdot \dfrac{25}{20}\left[-2P_2\left(\dfrac{3}{5}\right)+3P_2\left(\dfrac{6}{5}\right)+2P_2\left(\dfrac{9}{5}\right)-3P_2\left(\dfrac{12}{5}\right)\right]=1.\]
As the width of the cusp $3/10$ is $1$ by Lemma \ref{lem23}, we see that $\mbox{ord}_{3/10} l(2\tau) = 1\cdot 1=1$, so $l(2\tau)$ has a simple zero at $3/10$. Similarly, $l(2\tau)$ has another simple zero at $1/5$, and two simple poles at $1/10$ and $3/5$. We deduce that the total degree of poles of both $l(\tau)$ and $l(2\tau)$ are $2$, so by Proposition \ref{prop43}, there is a polynomial
\[L_2(X,Y) = \sum_{0\leq i, j\leq 2}C_{i,j}X^iY^j\in \mathbb{C}[X,Y]\]
such that $L_2(l(\tau),l(2\tau))=0$. In addition, we get that $C_{2,2}\neq 0, C_{0,0}\neq 0$ and $C_{2,0}=C_{2,1}=C_{0,1}=C_{0,2}=0$. Using the $q$-expansion (\ref{eq31})
of $l(\tau)$, we may take $C_{2,2}=1$. Hence, we arrive at $L_2(X,Y)= 1+X-2XY-XY^2+X^2Y^2$ and the desired identity follows.
\end{proof}

\begin{theorem}[Modular equation of level five]\label{thm46}
We have
\begin{align*}
(l(\tau)^5-l(5\tau))(l(5\tau)^2-l(5\tau)-1)^2=-5l(\tau)l(5\tau)P(l(\tau),l(5\tau)),
\end{align*}
where $P(X,Y) := (Y^2-Y-1)YX^3+(3Y^2+2Y-3)X^2-(3Y^2+2Y-3)YX+Y^2-Y-1$.
\end{theorem}

\begin{proof}
Appealing to Lemmas \ref{lem22} and \ref{lem42}, we focus on the congruence subgroup $\Gamma_1(10)\cap \Gamma_0(50)$ and the following cusps:
\[1/5, 2/5, 3/5, 4/5, 1/15, 2/15, 4/15, 8/15, 1/25,3/25.\]
Applying Lemmas \ref{lem23} and \ref{lem29}, we see that $l(\tau)$ has five simple poles at $1/5, 4/5, 1/15, 4/15$ and $1/25$, and five simple zeros at $2/5, 3/5, 2/15, 8/15$ and $3/25$. We now compute the order of $l(5\tau)$ at the cusp $1/25$ (which is a pole) as follows. Observe that
\[l\left(5\begin{bmatrix}
1 & 0\\25 & 1
\end{bmatrix}\tau\right)=l\left(\begin{bmatrix}
5 & 0\\25 & 1
\end{bmatrix}\tau\right)=l\left(\begin{bmatrix}
1 & 0\\5 & 1
\end{bmatrix}(5\tau)\right),\]
so by Lemma \ref{lem29}, the $q$-expansion of $l(5[\begin{smallmatrix}
1 & 0\\ 25 &1
\end{smallmatrix}]\tau)$ starts with $\varepsilon q^{\delta}$ with $|\varepsilon|=1$ and 
\[\delta = 5\cdot \dfrac{25}{20}\left[-2P_2\left(\dfrac{1}{5}\right)+3P_2\left(\dfrac{2}{5}\right)+2P_2\left(\dfrac{3}{5}\right)-3P_2\left(\dfrac{4}{5}\right)\right]=-\frac{5}{2}.\]
Since the width of the cusp $1/25$ is $2$ by Lemma \ref{lem23}, we see that $\mbox{ord}_{1/25} l(5\tau) = 2\cdot -5/2=-5$. Similarly, we get that $\mbox{ord}_{3/25} l(5\tau) = 5$, so the total degree of poles of both $l(\tau)$ and $l(5\tau)$ are $5$. Thus, by Proposition \ref{prop43}, there is a polynomial 
\[L_5(X,Y) = \sum_{0\leq i, j\leq 5}C_{i,j}X^iY^j\in \mathbb{C}[X,Y]\]
such that $L_5(l(\tau),l(5\tau))=0$. In addition, we get that $C_{5,0}\neq 0$ and $C_{0,1}\neq 0$. Switching the roles of $l(\tau)$ and $l(5\tau)$ yields $C_{0,5}\neq 0, C_{j,5}=0$ for $j\in \{1,2,3,4,5\}$ and $C_{k,0}=0$ for $k\in \{0,1,2,3,4\}$. Using the $q$-expansion (\ref{eq31}) of $l(\tau)$ leads us to 
\[L_5(X,Y) = (X^5-Y)(Y^2-Y-1)^2+5XYP(X,Y)\]
and the desired identity follows.
\end{proof}

\begin{theorem}\label{thm47}
Let $0=L_p(X,Y) = \sum_{0\leq i, j\leq p+1} C_{i,j} X^iY^j\in \mathbb{Q}[X,Y]$ be the modular equation for $l(\tau)$ of odd prime level $p\neq 5$.
\begin{enumerate}
\item[(1)] If $p\equiv \pm 1\pmod{10}$, then $C_{p+1,0}\neq 0, C_{0,p+1}\neq 0$, and $C_{i,p+1}=C_{p+1,i}=0$ for $i\in \{1,\ldots, p+1\}$ and $C_{0,j}=C_{j,0}=0$ for $j\in \{0,\ldots,p\}$.
\item[(2)] If $p\equiv \pm 3\pmod{10}$, then $C_{p+1,p}\neq 0, C_{0,1}\neq 0, C_{1,p+1}\neq 0, C_{p,0}\neq 0$, and $C_{p+1,i}=C_{i,0}=0$ for $i\in \{0,\ldots, p-1,p+1\}$ and $C_{0,j}=C_{j,p+1}=0$ for $j\in \{0,2,\ldots,p+1\}$.
\end{enumerate}
\end{theorem}

\begin{proof}
We adapt the notations of Lemma \ref{lem22}, and in view of Lemma \ref{lem42}, we need to determine the sets $S_c$ and $A_c$ for $c\in \{5,5p\}$. We set $I := \{0,1,\ldots,p-1\}$ and choose a unique integer $k_0\in I$ such that $p\mid 1+10k_0$. Then 
\begin{align*}
\Delta = \{\pm(1+10k)\in (\mathbb{Z}/10p\mathbb{Z})^\times : k\in I\setminus\{k_0\}\}.
\end{align*}
For each $a\in  I\setminus\{k_0\}$, choose $b\in I\setminus\{k_0\}$ such that $p\mid 1+5(a+b)$. Then $\pi_{2p}(1+10a)=\pi_{2p}(-1-10b)$, so $\pi_{2p}(\Delta)$ has exactly $p-1$ distinct elements modulo $2p$. We see that $\pi_{2p}(\Delta)=(\mathbb{Z}/2p\mathbb{Z})^\times$ and thus, $S_5=\{1\}$. On the other hand, we choose $d\in I\setminus\{k_0\}$ such that $p\mid 1+5d$. Then $\pi_{2p}(-1-10d)=1$, so $\Delta\cap \ker\pi_{2p}=\{1,-1-10d\}$ and $\pi_5(\Delta\cap \ker\pi_{2p})=\{\pm 1\}$. Thus, we take $A_5=\{1,2\}$. We next observe that $\pi_2(\Delta) = \{1\}=(\mathbb{Z}/2\mathbb{Z})^\times$, so we have $S_{5p}=\{1\}$ and  $\pi_{5p}(\Delta\cap \ker\pi_2)=\pi_{5p}(\Delta)\subset (\mathbb{Z}/5p\mathbb{Z})^\times$. We claim that $\pi_{5p}(\Delta)$ has $2(p-1)$ distinct elements modulo $5p$; otherwise, we would have $1+10r\equiv \pm(1+10s)\pmod{5p}$ for some distinct $r, s\in I\setminus\{k_0\}$, and all cases lead to a contradiction. As $(\mathbb{Z}/5p\mathbb{Z})^\times$ has $4(p-1)$ elements, we take $A_{5p}=\{1,2\}$.

We see that $l(\tau)$ and $l(p\tau)$ has poles and zeros at the following cusps of $\Gamma_1(10)\cap \Gamma_0(10p)$: $1/5, 2/5, 1/5p$ and $2/5p$. Applying Lemma \ref{lem29}, we get $\mbox{ord}_{1/5} l(\tau)=-p=-\mbox{ord}_{2/5} l(\tau)$ and $\mbox{ord}_{1/5p} l(\tau)=-1=-\mbox{ord}_{2/5p} l(\tau)$. We next compute $\mbox{ord}_{1/5p} l(p\tau)$ and  $\mbox{ord}_{2/5p} l(p\tau)$. We note that
\[l\left(p\begin{bmatrix}
1 & 0\\5p & 1
\end{bmatrix}\tau\right)=l\left(\begin{bmatrix}
p & 0\\5p & 1
\end{bmatrix}\tau\right)=l\left(\begin{bmatrix}
1 & 0\\5 & 1
\end{bmatrix}(p\tau)\right),\]
so by Lemma \ref{lem29}, the first term of the $q$-expansion of $l(p[\begin{smallmatrix}
1 & 0\\ 5p &1
\end{smallmatrix}]\tau)$ is $\varepsilon q^{\delta}$ with $|\varepsilon|=1$ and 
\[\delta = p\cdot \dfrac{25}{20}\left[-2P_2\left(\dfrac{1}{5}\right)+3P_2\left(\dfrac{2}{5}\right)+2P_2\left(\dfrac{3}{5}\right)-3P_2\left(\dfrac{4}{5}\right)\right]=-\frac{p}{2}.\]
Since the width of the cusp $1/5p$ is $2$ by Lemma \ref{lem23}, we get $\mbox{ord}_{1/5p}l(p\tau) = 2\cdot-p/2=-p$. On the other hand, we choose integers $u$ and $v$ such that $2v-5pu=1$. We have 
\[l\left(p\begin{bmatrix}
2 & u\\5p & v
\end{bmatrix}\tau\right)=l\left(\begin{bmatrix}
2p & up\\5p & v
\end{bmatrix}\tau\right)=l\left(\begin{bmatrix}
2 & up\\5 & v
\end{bmatrix}(p\tau)\right),\]
so Lemma \ref{lem29} shows that the $q$-expansion of $l(p[\begin{smallmatrix}
2 & u\\ 5p &v
\end{smallmatrix}]\tau)$ starts with $\varepsilon q^{\delta}$ with $|\varepsilon|=1$ and
\[\delta = p\cdot \dfrac{25}{20}\left[-2P_2\left(\dfrac{2}{5}\right)+3P_2\left(\dfrac{4}{5}\right)+2P_2\left(\dfrac{6}{5}\right)-3P_2\left(\dfrac{8}{5}\right)\right]=\frac{p}{2}.\]
As the cusp $2/5p$ has width $2$ by Lemma \ref{lem23}, we obtain $\mbox{ord}_{2/5p}l(p\tau) = 2\cdot p/2=p$.

We now compute $\mbox{ord}_{1/5} l(p\tau)$ and  $\mbox{ord}_{2/5} l(p\tau)$. We choose an integer $g$ with $pg\equiv 1\pmod{5}$, so that
\[l\left(p\begin{bmatrix}
1 & 0\\5 & 1
\end{bmatrix}\tau\right)=l\left(\begin{bmatrix}
p & 0\\5 & 1
\end{bmatrix}\tau\right)=l\left(\begin{bmatrix}
p & (pg-1)/5\\5 & g
\end{bmatrix}\dfrac{\tau+(1-pg)/5}{p}\right).\]
By Lemma \ref{lem29}, we deduce that the $q$-expansion of $l(p[\begin{smallmatrix}
1 & 0\\ 5 &1
\end{smallmatrix}]\tau)$ starts with $\varepsilon q^{\delta_p}$ with $|\varepsilon|=1$ and
\[\delta_p := \dfrac{1}{p}\cdot \dfrac{25}{20}\left[-2P_2\left(\dfrac{p}{5}\right)+3P_2\left(\dfrac{2p}{5}\right)+2P_2\left(\dfrac{3p}{5}\right)-3P_2\left(\dfrac{4p}{5}\right)\right].\]
Since the width of $1/5$ is $2p$ by Lemma \ref{lem23}, we see that
\[\mbox{ord}_{1/5}l(p\tau) = 2p\cdot \delta_p=\begin{cases}
-1 &\text{ if }p\equiv \pm 1\pmod{10},\\
1&\text{ if }p\equiv \pm 3\pmod{10}.
\end{cases}\]
We next take an integer $h$ with $2ph\equiv 1\pmod{5}$, so that
\[l\left(p\begin{bmatrix}
2 & 1\\5 & 3
\end{bmatrix}\tau\right)=l\left(\begin{bmatrix}
2p & p\\5 & 3
\end{bmatrix}\tau\right)=l\left(\begin{bmatrix}
2p & (2ph-1)/5\\5 & h
\end{bmatrix}\dfrac{\tau+(3-ph)/5}{p}\right).\]
Using Lemma \ref{lem29}, the first term of the $q$-expansion of $l(p[\begin{smallmatrix}
2 & 3\\ 5 &1
\end{smallmatrix}]\tau)$ is $\varepsilon q^{\delta_p}$ with $|\varepsilon|=1$ and
\[\delta_p := \dfrac{1}{p}\cdot \dfrac{25}{20}\left[-2P_2\left(\dfrac{2p}{5}\right)+3P_2\left(\dfrac{4p}{5}\right)+2P_2\left(\dfrac{6p}{5}\right)-3P_2\left(\dfrac{8p}{5}\right)\right].\]
Since the width of $2/5$ is $2p$ by Lemma \ref{lem23}, we arrive at
\[\mbox{ord}_{2/5}l(p\tau) = 2p\cdot \delta_p=\begin{cases}
1 &\text{ if }p\equiv \pm 1\pmod{10},\\
-1&\text{ if }p\equiv \pm 3\pmod{10}.
\end{cases}\]
Using the notations of Proposition \ref{prop43} with $f_1(\tau)=l(\tau)$ and $f_2(\tau)=l(p\tau)$, we obtain $S_{1,\infty}=\{1/5,1/5p\}, S_{1,0}=\{2/5,2/5p\}$ and 
\begin{align*}
S_{2,\infty} &= \begin{cases}
\{1/5, 1/5p\}&\text{ if }p\equiv \pm 1\pmod{10},\\
\{2/5, 1/5p\}&\text{ if }p\equiv \pm 3\pmod{10},
\end{cases}\\
S_{2,0} &= \begin{cases}
\{2/5, 2/5p\}&\text{ if }p\equiv \pm 1\pmod{10},\\
\{1/5, 2/5p\}&\text{ if }p\equiv \pm 3\pmod{10},
\end{cases}
\end{align*}
so that $S_{1,\infty}\cup S_{1,0}=S_{2,\infty}\cup S_{2,0}$. In addition, we deduce that the total number of poles of both $l(\tau)$ and $l(p\tau)$ is $p+1$, so there is a polynomial 
\[L_p(X,Y) = \sum_{0\leq i, j\leq p+1} C_{i,j} X^iY^j\in \mathbb{Q}[X,Y]\]
with $L_p(l(\tau),l(p\tau))=0$. We find that if $p\equiv \pm 1\pmod{10}$, then $C_{p+1,0}\neq 0$ and $C_{0,p+1}\neq 0$, so that $C_{p+1,i}=0$ for $i\in \{1,\ldots, p+1\}$ and $C_{0,j}=0$ for $j\in \{0,\ldots,p\}$. Switching the roles of $l(\tau)$ and $l(p\tau)$ yields $C_{i,p+1}=0$ for $i\in \{1,\ldots, p+1\}$ and $C_{j,0}=0$ for $j\in \{0,\ldots,p\}$. On the other hand, if $p\equiv \pm 3\pmod{10}$, then $C_{p+1,p}\neq 0$ and $C_{0,1}\neq 0$, implying that  $C_{p+1,i}=0$ for $i\in \{0,\ldots, p-1,p+1\}$ and $C_{0,j}=0$ for $j\in \{0,2,\ldots,p+1\}$. Swapping the roles of $l(\tau)$ and $l(p\tau)$ gives $C_{1,p+1}\neq 0, C_{p,0}\neq 0, C_{i,0}=0$ for $i\in \{0,\ldots, p-1,p+1\}$ and $C_{j,p+1}=0$ for $j\in \{0,2,\ldots,p+1\}$.
\end{proof}

We give the modular equations $L_n(X,Y)$ for $l(\tau)$ of level $n$ with $n\in \{2,4,5,6,13\}$ as shown in Table \ref{tab:tbl2}. On the other hand, we verify via \textit{Mathematica} that the modular equations for $l(\tau)$ of prime levels $p\in \{3,7,11,19\}$ are the same as those for $k(\tau)$, which appeared as a supplementary material (Appendix A) to \cite{leepark}, so we will not present them here and instead remark that they satisfy the Kronecker congruence depending on the values of $p$ modulo $10$.

\begin{table}[h]
\begin{center}
{\tiny\begin{tabular}{@{}ll@{}}
	\toprule
	$n$ & $L_n(X,Y)$ \\
	\midrule
	$2$ & $1 + X - 2 X Y - X Y^2 + X^2 Y^2$ \\
	\midrule
	$4$ & $\begin{aligned}
		&-X - X^2 + X^3 + X^4 + 4 X Y + 4 X^2 Y - 8 X^3 Y + 2 X Y^2 + 8 X^2 Y^2 - 2 X^3 Y^2- 8 X Y^3- 4 X^2 Y^3 + 4 X^3 Y^3\\
		&+ Y^4 - X Y^4 - X^2 Y^4 + X^3 Y^4
	\end{aligned}$\\
	\midrule
	$5$ & $\begin{aligned}
		&(X^5-Y)(Y^2-Y-1)^2+5 X Y (1 - 3 X^2 + Y + 3 X Y + 2 X^2 Y - X^3 Y - Y^2 - 2 X Y^2+ 3 X^2 Y^2 - X^3 Y^2\\
		&- 3 X Y^3 + X^3 Y^3)
	\end{aligned}$\\
	\midrule
	$6$ & $\begin{aligned}
		&-X - 2 X^2 + 2 X^4 + X^5 + 6 X Y - 12 X^3 Y + 8 X^4 Y - 2 X^5 Y - 3 X Y^2 + 33 X^2 Y^2+ 27 X^3 Y^2 - 53 X^4 Y^2\\
		&- 13 X^5 Y^2 + 3 X^6 Y^2 + 5 X^7 Y^2 + X^8 Y^2 - 22 X Y^3 + 24 X^2 Y^3+ 96 X^3 Y^3 - 40 X^4 Y^3 - 78 X^5 Y^3 + 24 X^6 Y^3\\
		&+ 12 X^7 Y^3 + 9 X Y^4 - 30 X^2 Y^4- 39 X^3 Y^4 + 92 X^4 Y^4 + 39 X^5 Y^4 - 30 X^6 Y^4 - 9 X^7 Y^4 + 12 X Y^5 - 24 X^2 Y^5\\
		&- 78 X^3 Y^5 + 40 X^4 Y^5 + 96 X^5 Y^5 - 24 X^6 Y^5 - 22 X^7 Y^5 + Y^6 - 5 X Y^6 + 3 X^2 Y^6+ 13 X^3 Y^6- 53 X^4 Y^6\\
		&- 27 X^5 Y^6 + 33 X^6 Y^6 + 3 X^7 Y^6 - 2 X^3 Y^7 - 8 X^4 Y^7- 12 X^5 Y^7 + 6 X^7 Y^7 - X^3 Y^8 + 2 X^4 Y^8 - 2 X^6 Y^8\\
		&+ X^7 Y^8
	\end{aligned}$\\
	\midrule
	$13$ & $\begin{aligned}
		&(X^{13}-Y)(XY^{13}+1)+13 X Y(1 - 4 X + X^2 + 26 X^3 - 33 X^4 - 57 X^5 + 88 X^6 + 53 X^7- 77 X^8 - 18 X^9\\
		&+ 21 X^{10} + X^{11} - X^{12} + Y + 6 X Y - 88 X^2 Y + 76 X^3 Y + 580 X^4 Y- 214 X^5 Y - 990 X^6 Y + 226 X^7 Y\\
		&+ 608 X^8 Y - 100 X^9 Y - 110 X^{10} Y + 6 X^{11} Y+ 4 X^{12} Y - 21 Y^2 + 110 X Y^2 + 178 X^2 Y^2 - 1147 X^3 Y^2\\
		&- 629 X^4 Y^2 + 3003 X^5 Y^2+396 X^6 Y^2 - 2915 X^7 Y^2 + 255 X^8 Y^2 + 1053 X^9 Y^2 - 178 X^{10} Y^2 - 88 X^{11} Y^2\\
		&- X^{12} Y^2 - 18 Y^3 - 100 X Y^3 + 1053 X^2 Y^3 + 114 X^3 Y^3 - 6447 X^4 Y^3 - 586 X^5 Y^3+11132 X^6 Y^3 + 382 X^7 Y^3\\
		&- 6808 X^8 Y^3 + 114 X^9 Y^3 + 1147 X^{10} Y^3 + 76 X^{11} Y^3- 26 X^{12} Y^3 + 77 Y^4 - 608 X Y^4 - 255 X^2 Y^4 + 6808 X^3 Y^4\\
		&+ 1291 X^4 Y^4 - 18366 X^5 Y^4-484 X^6 Y^4 + 18077 X^7 Y^4 - 1291 X^8 Y^4 - 6447 X^9 Y^4 + 629 X^{10} Y^4 + 580 X^{11} Y^4 \\
		&+ 33 X^{12} Y^4 + 53 Y^5 + 226 X Y^5 - 2915 X^2 Y^5 + 382 X^3 Y^5 + 18077 X^4 Y^5 + 96 X^5 Y^5-30602 X^6 Y^5 + 96 X^7 Y^5\\
		&+ 18366 X^8 Y^5 - 586 X^9 Y^5 - 3003 X^{10} Y^5 - 214 X^{11} Y^5+ 57 X^{12} Y^5 - 88 Y^6 + 990 X Y^6 -396 X^2 Y^6\\
		&- 11132 X^3 Y^6 + 484 X^4 Y^6 + 30602 X^5 Y^6-30602 X^7 Y^6 + 484 X^8 Y^6 + 11132 X^9 Y^6 - 396 X^{10} Y^6 - 990 X^{11} Y^6\\ 
		&- 88 X^{12} Y^6- 57 Y^7 - 214 X Y^7 + 3003 X^2 Y^7 - 586 X^3 Y^7 - 18366 X^4 Y^7 + 96 X^5 Y^7 + 30602 X^6 Y^7+ 96 X^7 Y^7\\
		&- 18077 X^8 Y^7 + 382 X^9 Y^7 + 2915 X^{10} Y^7 + 226 X^{11} Y^7 - 53 X^{12} Y^7+ 33 Y^8 - 580 X Y^8 + 629 X^2 Y^8 + 6447 X^3 Y^8\\
		&- 1291 X^4 Y^8 - 18077 X^5 Y^8 - 484 X^6 Y^8+18366 X^7 Y^8 + 1291 X^8 Y^8 - 6808 X^9 Y^8 - 255 X^{10} Y^8 + 608 X^{11} Y^8\\
		&+ 77 X^{12} Y^8+ 26 Y^9 + 76 X Y^9 - 1147 X^2 Y^9 + 114 X^3 Y^9 + 6808 X^4 Y^9 + 382 X^5 Y^9 - 11132 X^6 Y^9-586 X^7 Y^9\\
		&+ 6447 X^8 Y^9 + 114 X^9 Y^9 - 1053 X^{10} Y^9 -100 X^{11} Y^9 + 18 X^{12} Y^9 - Y^{10}+ 88 X Y^{10} - 178 X^2 Y^{10} - 1053 X^3 Y^{10}\\ 
		&+ 255 X^4 Y^{10} + 2915 X^5 Y^{10} + 396 X^6 Y^{10}-3003 X^7 Y^{10} - 629 X^8 Y^{10} + 1147 X^9 Y^{10} + 178 X^{10} Y^{10} -110 X^{11} Y^{10}\\
		&- 21 X^{12} Y^{10}- 4 Y^{11} + 6 X Y^{11} + 110 X^2 Y^{11} - 100 X^3 Y^{11} - 608 X^4 Y^{11} + 226 X^5 Y^{11} + 990 X^6 Y^{11}-214 X^7 Y^{11}\\
		&- 580 X^8 Y^{11} + 76 X^9 Y^{11} + 88 X^{10} Y^{11} +6 X^{11} Y^{11} - X^{12} Y^{11} - Y^{12}- X Y^{12} + 21 X^2 Y^{12} + 18 X^3 Y^{12} - 77 X^4 Y^{12}\\
		&- 53 X^5 Y^{12} + 88 X^6 Y^{12} +57 X^7 Y^{12}- 33 X^8 Y^{12} - 26 X^9 Y^{12} + X^{10} Y^{12} + 4 X^{11} Y^{12} + X^{12} Y^{12})
	\end{aligned}$\\
	\bottomrule
\end{tabular}}
\end{center}
\caption{Modular equations $L_n(X,Y)=0$ for $l(\tau)$ of level $n\in\{2,4,5,6,13\}$}\label{tab:tbl2}%
\end{table}

We now consider modular equations for $l(\tau)$ of levels $n$ coprime to $10$. For any integer $a$ coprime to $10$, we choose a matrix $\sigma_a\in\mbox{SL}_2(\mathbb{Z})$ such that $\sigma_a\equiv [\begin{smallmatrix}
a^{-1} &0\\ 0 &a
\end{smallmatrix}]\pmod{10}$. We take
\begin{equation*}
\sigma_{\pm 1} =\pm\begin{bmatrix}
1 & 0\\ 0 & 1
\end{bmatrix}\quad\text{ and }\quad\sigma_{\pm 3}=\pm\begin{bmatrix}
-3 & -10\\10 & 33
\end{bmatrix}.
\end{equation*}
We know that $l(\sigma_{\pm 1}\tau) = l(\tau)$ and $l(\sigma_{\pm 3}\tau)=-1/l(\tau)$ by Lemma \ref{lem27}. Moreover, by \cite[Prop. 3.36]{shim}, we have the disjoint union 
\begin{align*}
\Gamma_1(10)\begin{bmatrix}
1 & 0\\0 & n
\end{bmatrix}\Gamma_1(10) = \bigsqcup_{0<a\mid n}\bigsqcup_{\substack{0\leq b < n/a\\\gcd(a,b,n/a)=1}}\Gamma_1(10)\sigma_a\begin{bmatrix}
a & b\\0 & \frac{n}{a}
\end{bmatrix}
\end{align*}
and $\Gamma_1(10)\backslash \Gamma_1(10)[\begin{smallmatrix}
1 & 0\\0 & n
\end{smallmatrix}]\Gamma_1(10)$ has $\psi(n):=n\prod_{p\mid n}(1+1/p)$ elements. We now define the polynomial
\begin{align*}
\Phi_n(X,\tau):= \prod_{0<a\mid n}\prod_{\substack{0\leq b < n/a\\\gcd(a,b,n/a)=1}}(X-(l\circ\alpha_{a,b})(\tau))
\end{align*}
where $\alpha_{a,b}:=\sigma_a[\begin{smallmatrix}
a & b\\0 & n/a
\end{smallmatrix}]$. Then the coefficients of $\Phi_n(X,\tau)$ are elementary symmetric functions of $l\circ\alpha_{a,b}$, so these are invariant under the action of $\Gamma_1(10)$ and thus are in $A_0(\Gamma_1(10))=\mathbb{C}(l(\tau))$. Thus, we may write $\Phi_n(X,\tau)$ as $\Phi_n(X,l(\tau))\in\mathbb{C}(l(\tau))[X]$. 

Let $S_{n,\infty}$ (resp. $S_{n,0}$) be the set of inequivalent cusps of $\Gamma_1(10)\cap \Gamma_0(10n)$ for which $l(n\tau)$ has a pole (resp. a zero). We define 
\begin{align*}
r_n := -\sum_{r\in S_{1,\infty}\cap S_{n,0}}\mbox{ord}_r l(\tau),
\end{align*}
with the convention that $r_n=0$ if $S_{1,\infty}\cap S_{n,0}$ is empty, and set 
\begin{align*}
L_n(X,l(\tau)) = l(\tau)^{r_n} \Phi_n(X,l(\tau)).
\end{align*}
The following result gives the properties of the polynomial $L_n(X,l(\tau))$ when $\gcd(n,10)=1$. We omit the proof as it is similar to that of \cite[Thm. 1.1]{leepark}. We note that $L_n(l(\frac{\tau}{n}),l(\tau))=0$, which can be seen as the modular equation for $l(\tau)$ of level $n$.

\begin{theorem}\label{thm48}
Let $n$ be a positive integer coprime to $10$ and let $L_n(X,Y)$ be the polynomial as above. Then
\begin{enumerate}
\item[(1)] $L_n(X,Y)\in \mathbb{Z}[X,Y]$ and $\deg_X L_n(X,Y)=\psi(n)$.
\item[(2)] $L_n(X,Y)$ is irreducible both as a polynomial in $X$ over $\mathbb{C}(Y)$ and as a polynomial in $Y$ over $\mathbb{C}(X)$.
\item[(3)] We have 
\begin{align*}
L_n(X,Y)=\begin{cases}
	L_n(Y,X) &\text{ if }n\equiv \pm 1\pmod{10},\\
	Y^{\psi(n)}L_n(-Y^{-1},X)&\text{ if }n\equiv \pm 3\pmod{10}.
\end{cases}
\end{align*}
\item[(4)] (Kronecker congruence) If $p\neq 5$ is an odd prime, then 
\begin{align*}
L_p(X,Y)\equiv\begin{cases}
	(X^p-Y)(X-Y^p)\pmod{p\mathbb{Z}[X,Y]} &\text{ if }p\equiv \pm 1\pmod{10},\\
	(X^p-Y)(XY^p+1)\pmod{p\mathbb{Z}[X,Y]}&\text{ if }p\equiv \pm 3\pmod{10}.
\end{cases}
\end{align*}
\end{enumerate}
\end{theorem}

\section{Arithmetic of \texorpdfstring{$l(\tau)$}{l(tau)} and evaluation of \texorpdfstring{$r(\tau)^ar(2\tau)^b$}{r(t)ar(t)b}}\label{sec5}

We study in this section some arithmetic properties of $l(\tau)$ at imaginary quadratic points $\tau\in\mathbb{H}$, which will be useful in evaluating the functions of the form $r(\tau)^ar(2\tau)^b$ satisfying the conditions of Theorem \ref{thm31}. We first state the following result of Cho and Koo \cite{chokoo} about generating ray class fields by a singular value of a modular function on some congruence subgroup between $\Gamma(N)$ and $\Gamma_1(N)$.

\begin{lemma}[\cite{chokoo}]\label{lem51}
Let $K$ be an imaginary quadratic field with discriminant $d_K$ and $\tau\in K\cap\mathbb{H}$ be a root of a primitive equation $aX^2+bX+c=0$ such that $b^2-4ac=d_K$ and $a,b,c\in\mathbb{Z}$. Let $\Gamma'$ be a congruence subgroup such that $\Gamma(N)\subset \Gamma'\subset \Gamma_1(N)$. Suppose that $\gcd(a,N)=1$. Then the field generated over $K$ by all values of $h(\tau)$, where $h\in A_0(\Gamma')_\mathbb{Q}$ is defined and finite at $\tau$, is the ray class field $K_{(N)}$ modulo $N$ over $K$.
\end{lemma}

We now apply Lemma \ref{lem51} to prove the following.

\begin{theorem}\label{thm52}
Let $K$ be an imaginary quadratic field with discriminant $d_K$ and $\tau\in K\cap\mathbb{H}$ be a root of the primitive polynomial $aX^2+bX+c\in\mathbb{Z}[X]$ such that $b^2-4ac=d_K$ and $\gcd(a,10)=1$. Then $K(l(\tau))$ is the ray class field modulo $10$ over $K$. 
\end{theorem}

\begin{proof}
In view of Theorem \ref{thm31}, we have that $A_0(\Gamma_1(10))_{\mathbb{Q}}=\mathbb{Q}(l(\tau))$ since $l(\tau)$ has a rational $q$-expansion. Suppose $\tau\in K\cap\mathbb{H}$ is a root of the primitive polynomial $aX^2+bX+c\in\mathbb{Z}[X]$ such that $b^2-4ac=d_K$ and $\gcd(a,10)=1$. Then $l(\tau)$ is well-defined and finite, so the desired conclusion follows from Lemma \ref{lem51}.
\end{proof}

The following result shows that the modular $j$-invariant defined by 
\begin{align*}
j(\tau) = \eta^{-24}(\tau)\left(1+240\sum_{n=1}^\infty\dfrac{n^3q^n}{1-q^n}\right)^3
\end{align*}
can be written as a rational function of $l(\tau)$, which implies that $l(\tau)$ is an algebraic unit for imaginary quadratic points $\tau\in\mathbb{H}$. In general, we infer from \cite[Prop. 3]{geehons} that the functions $r(\tau)^ar(2\tau)^b$ are algebraic units for all integers $a$ and $b$ and for such points $\tau\in\mathbb{H}$. 

\begin{theorem}\label{thm53}
Let $K$ be an imaginary quadratic field and $\tau\in K\cap\mathbb{H}$. Then
\begin{align*}
j(\tau) = \dfrac{P(l)^3}{l^2(l^2-1)(l^2+4l-1)^5(l^2-l-1)^{10}},
\end{align*}
where $l:=l(\tau)$ and 
\begin{align*}
P(X) := 1 &- 4 X + 234 X^2 - 460 X^3 + 495 X^4 + 456 X^5 - 1444 X^6 \\
&- 456 X^7 + 495 X^8 + 460 X^9 + 234 X^{10} + 4 X^{11} + X^{12}.
\end{align*}
\end{theorem}

\begin{proof}
We first remark that the function $f(\tau) := \eta^6(\tau)\eta^{-6}(5\tau)$ generates the field of all modular functions on $\Gamma_0(5)$ by Lemmas \ref{lem25} and \ref{lem26}, and satisfies the identity \cite[Thm. 5.26]{cooper}
\begin{align}
j(\tau) = \dfrac{(f^2(\tau)+250f(\tau)+3125)^3}{f^5(\tau)}.\label{eq51}
\end{align}As a modular function on $\Gamma_0(10)$, we deduce from Lemma \ref{lem26} that $f(\tau)$ has a double pole at $1/5$ and a simple pole at $\infty$. Since $A_0(\Gamma_0(10))=\mathbb{C}(f(\tau),l(\tau)-1/l(\tau))$ by Theorem \ref{thm32}, we infer from Proposition \ref{prop43} that there is a polynomial
\begin{align*}
F(X,Y) = \sum_{\substack{0\leq i\leq 1\\0\leq j\leq 3}}C_{i,j}X^iY^j\in\mathbb{C}[X,Y]
\end{align*}
such that $F(f(\tau),l(\tau)-1/l(\tau))=0$. Using (\ref{eq31}) and the $q$-expansion of $f(\tau)$ leads to $F(X,Y)=4 - 7 Y - X Y + 2 Y^2 + Y^3$ and 
\begin{align}
f(\tau) = \dfrac{(l^2-l-1)^2(l^2+4l-1)}{l^2(l^2-1)},\label{eq52}
\end{align}
where $l:=l(\tau)$. Combining (\ref{eq51}) and (\ref{eq52}) gives the desired identity.
\end{proof}

The next result concerns about the singular values of $l(\rho\tau)$ for positive rational numbers $\rho$, which can be written as radicals using the values of $l(\tau)$. We omit the proof as it is analogous to that of \cite[Thm. 1.4]{leepark}.

\begin{theorem}\label{thm54}
For any positive $\rho\in\mathbb{Q}$, if $l(\tau)$ can be expressed in terms of radicals, then $l(\rho\tau)$ can also be expressed in terms of radicals.
\end{theorem}

We now discuss a method of evaluating singular values of $l(\theta)$, where $\theta\in K\cap \mathbb{H}$ is a generator for the ring of integers $\mathcal{O}_K$ of an imaginary quadratic field $K$ with discriminant $d_K\equiv 0\pmod{4}$, by finding its minimal polynomial over $K$, which is the class polynomial of the ray class field $K_{(10)}$ modulo $10$ over $K$, using the Shimura reciprocity law. For more details, we refer the reader to \cite{chokoo, gee, guad}.

We denote by $\mathcal{F}_N$ the field modular functions of level $N$ whose coefficients in their $q$-expansions lie in the cyclotomic field $\mathbb{Q}(\zeta_N)$. 
We know from the main theorem of complex multiplication that if $h\in\mathcal{F}_N$ is a modular function with $h(\theta)$ finite, then $h(\theta)$ generates the ray class field $K_{(N)}$ over $K$.
Consider the map $g_\theta : (\mathcal{O}_K/N\mathcal{O}_K)^\times\rightarrow \mbox{GL}_2(\mathbb{Z}/N\mathbb{Z})$ given by
\begin{equation*}
g_\theta(s\theta+t) = \begin{bmatrix}
t - Bs & -Cs\\ s & t
\end{bmatrix}\in \mbox{GL}_2(\mathbb{Z}/N\mathbb{Z}),
\end{equation*}
where $s\theta+t\in (\mathcal{O}_K/N\mathcal{O}_K)^\times$ and $X^2+BX+C\in\mathbb{Z}[X]$ is the minimal polynomial of $\theta$ over $\mathbb{Q}$. Let $K(j(\theta))$ be the Hilbert class field of $K$ and denote by $W_{N,\theta}$ the image of $(\mathcal{O}_K/N\mathcal{O}_K)^\times$ under $g_\theta$. The Shimura reciprocity law states that there is a surjective homomorphism from $W_{N,\theta}$ to $\mbox{Gal}(K_{(N)}/K(j(\theta)))$ that sends $\alpha$ to the map $h(\theta)\mapsto h^{\alpha^{-1}}(\theta)$, whose kernel $T$ is the image $g_\theta(\mathcal{O}_K^\times)$ given by
\begin{equation*}
T = \begin{cases}
\left\{\pm\begin{bmatrix}
1 & 0\\ 0 & 1
\end{bmatrix}, \pm\begin{bmatrix}
0 &-1\\ 1 & 0
\end{bmatrix}\right\} & \text{ if } K=\mathbb{Q}(\sqrt{-1}),\\
\noalign{\vskip6pt}
\left\{\pm\begin{bmatrix}
1 & 0\\ 0 & 1
\end{bmatrix}, \pm\begin{bmatrix}
-1 &-1\\ 1 & 0
\end{bmatrix}, \pm\begin{bmatrix}
0 &-1\\ 1 & 1
\end{bmatrix}\right\} & \text{ if } K=\mathbb{Q}(\sqrt{-3}),\\
\noalign{\vskip6pt}
\left\{\pm\begin{bmatrix}
1 & 0\\ 0 & 1
\end{bmatrix}\right\} & \text{ otherwise.}
\end{cases}
\end{equation*}
To get the conjugates of $h(\theta)\in K_{(N)}$ over $K$, we now exploit the action of the form class group $C(d_K)$ on $h(\theta)$ as follows. Suppose $x:=[a,b,c]\in C(d_K)$ is a reduced primitive binary quadratic form and $\tau_x := (-b+\sqrt{d_K})/(2a)\in\mathbb{H}$. We choose a matrix $u_x\in\mbox{GL}_2(\mathbb{Z}/N\mathbb{Z})$ such that
\begin{equation*}
u_x \equiv\begin{cases}
\begin{bmatrix}
a & b/2\\ 0 & 1
\end{bmatrix} \bmod{p^{r_p}} & \text{ if } \gcd(p,a)=1,\\
\noalign{\vskip6pt}
\begin{bmatrix}
-b/2 & -c\\ 1 & 0
\end{bmatrix} \bmod{p^{r_p}} & \text{ if } p \mid a\text{ and } \gcd(p,c)=1,\\
\noalign{\vskip6pt}
\begin{bmatrix}
-b/2-a &-b/2-c\\ 1 & -1
\end{bmatrix} \bmod{p^{r_p}} & \text{ if } p \mid a\text{ and } p\mid c
\end{cases}
\end{equation*}
for each prime power $p^{r_p}$ dividing $N$. We have the action $h(\theta)^{x^{-1}} = h^{u_x}(\tau_x)$, and note that if $f$ is a modular function such that $f(\theta)\in K(j(\theta))$, then there is an isomorphism from $C(d_K)$ to $\mbox{Gal}(K(j(\theta))/K)$ that sends the class $x^{-1}$ to the map $f(\theta)\mapsto f^{u_x}(\tau_x)$. In view of the isomorphism 
\begin{equation*}
\mbox{Gal}(K_{(N)}/K)/\mbox{Gal}(K_{(N)}/K(j(\theta)))\simeq\mbox{Gal}(K(j(\theta))/K),
\end{equation*}
we obtain the following.

\begin{proposition}[\cite{chokoo}]\label{prop55}
Using the above notations, the set of all conjugates of $h(\theta)\in K_{(N)}$ over $K$ is given by
\begin{equation*}
\{ h^{\alpha\cdot u_x}(\tau_x) : \alpha \in W_{N,\theta}/T, x\in C(d_K)\}.
\end{equation*}
\end{proposition}

After obtaining all conjugates of $h(\theta)$, we now define the class polynomial of $K_{(N)}$ given by 
\begin{equation}\label{eq56}
F_N(X) := \prod_{\substack{\alpha \in W_{N,\theta}/T\\x\in C(d_K)}}(X-h^{\alpha\cdot u_x}(\tau_x))\in K[X].
\end{equation}

We now let $\theta := \frac{\sqrt{d_K}}{2}$; then $l(\theta)$ is a real number since $e^{2\pi i\sqrt{d_K}/2}\in\mathbb{R}$. Thus, 
\begin{align*}
0 = F_{10}(l(\theta))=\overline{F_{10}(l(\theta))}=\overline{F_{10}}(\overline{l(\theta)})=\overline{F_{10}}(l(\theta)),
\end{align*}
so that $F_{10}(X) \in (K\cap \mathbb{R})[X]  = \mathbb{Q}[X]$. We see from Theorem \ref{thm53} that $l(\theta)$ is an algebraic integer, which shows that $F_{10}(X) \in \mathbb{Z}[X]$. In this case, we can numerically obtain $F_{10}(X)$ using an accurate approximation of $l(\theta)$ via (\ref{eq31}). We illustrate an example of finding $F_{10}(X)$, which can be used to evaluate $l(\theta)$ and the functions $r(\theta)^ar(2\theta)^b$ using Corollary \ref{cor34}.

\begin{example}
Consider the field $K=\mathbb{Q}(i)$ with $d_K=-4$ and $\mathcal{O}_K=\mathbb{Z}[i]$, where $i=\sqrt{-1}$. We wish to find the class polynomial $F_{10}(X)$ with $F_{10}(l(i))=0$. We choose $(B,C) = (0,1)$ since $X^2+1$ is the minimal polynomial of $i$. Since $x:=[1,0,1]$ is the only reduced primitive binary quadratic form of discriminant $-4$, we have $\tau_x=i$ and $u_x = [\begin{smallmatrix} 1& 0\\0&1\end{smallmatrix}]\in \mbox{GL}_2(\mathbb{Z}/10\mathbb{Z})$. We compute
\begin{align*}
\begin{split}
W_{10,\theta}/T &= \left\lbrace
\begin{bmatrix}
	1 & 0\\0 & 1
\end{bmatrix},
\begin{bmatrix}
	3 & 0\\0 & 3
\end{bmatrix},
\begin{bmatrix}
	-1 & -4\\4 & -1
\end{bmatrix},
\begin{bmatrix}
	2 & -3\\3 & 2
\end{bmatrix},\right.\\
&\left.\quad\begin{bmatrix}
	5 & 2\\-2 & 5
\end{bmatrix},
\begin{bmatrix}
	-4 & -5\\5 & -4
\end{bmatrix},
\begin{bmatrix}
	-3 & 2\\-2 & -3
\end{bmatrix},
\begin{bmatrix}
	4 & 1\\-1 & 4
\end{bmatrix}\right\rbrace.
\end{split}
\end{align*}	
Using the $q$-expansion (\ref{eq31}) of $l(\tau)$, we get
\begin{align*}
F_{10}(X) &= (X-l(i))(X-l^{[\begin{smallmatrix} 3& -10\\10&-33\end{smallmatrix}]}(i))(X-l^{[\begin{smallmatrix} -1& -4\\2&7\end{smallmatrix}]}(i))(X-l^{[\begin{smallmatrix} 2& 7\\1&4\end{smallmatrix}]}(i))\\
&\times(X-l^{[\begin{smallmatrix} 5& 12\\2&5\end{smallmatrix}]}(i)(X-l^{[\begin{smallmatrix} 6& -5\\5&-4\end{smallmatrix}]}(i))(X-l^{[\begin{smallmatrix} 7& 2\\-4&-1\end{smallmatrix}]}(i))(X-l^{[\begin{smallmatrix} 4& 1\\7&2\end{smallmatrix}]}(i))\\
&\approx X^8+26X^7+62X^6+458X^5-130X^4-458X^3+62X^2-26X+1.
\end{align*}
As $l(i)\approx 1.00373486$, we deduce that
\begin{align*}
l(i) =\dfrac{1}{4}\left(-13-5\sqrt{5}+(15-\sqrt{5})\alpha+2\sqrt{5}\sqrt{25+13\sqrt{5}-(15+5\sqrt{5})\alpha}\right),
\end{align*}
where $\alpha:=\sqrt{2+\sqrt{5}}$. In view of Corollary \ref{cor34} and the approximations $r(2i)^3/r(i)\approx 0.00187091$ and $r(2i)/r(i)^7\approx 542.52907744$, we arrive at the exact values of $r(2i)^3/r(i)$ and $r(2i)/r(i)^7$ respectively given by
{\footnotesize \begin{align*}
& \dfrac{1}{4}\left(251+115\sqrt{5}-(115+59\sqrt{5})\alpha+2\sqrt{5}\sqrt{15725+7033\sqrt{5}-(7645+3415\sqrt{5})\alpha}\right),\\
&\dfrac{1}{4}\left(-6259-2815\sqrt{5}+(3145+1557\sqrt{5})\alpha+26\sqrt{5}\sqrt{50825+22733\sqrt{5}-(24615+11005\sqrt{5})\alpha}\right).
\end{align*}}%
On the other hand, using the modular equation for $l(\tau)$ of level two from Theorem \ref{thm45} and $l(i/2)\approx 1.08609902$, we obtain the
exact value of $l(i/2)$:
{\small \begin{align*}
\dfrac{1}{2}\left(-134-60\sqrt{5}-(65+29\sqrt{5})\alpha+\sqrt{10}\sqrt{7225+3231\sqrt{5}+(3510+1570\sqrt{5})\alpha}\right).
\end{align*}}%
We can then obtain the exact values of $r(i)^3/r(i/2)$ and $r(i)/r(i/2)^7$ using Corollary \ref{cor34} and the value of $l(i/2)$.
\end{example}

\section*{Acknowledgment}

The author would like to thank the anonymous referee for helpful suggestions that improved the quality of this paper.

\bibliography{prodrogram}

\providecommand{\bysame}{\leavevmode\hbox to3em{\hrulefill}\thinspace}
\providecommand{\MR}{\relax\ifhmode\unskip\space\fi MR }
\providecommand{\MRhref}[2]{%
  \href{http://www.ams.org/mathscinet-getitem?mr=#1}{#2}
}
\providecommand{\href}[2]{#2}
\begin{thebibliography}{10}

\bibitem{andber}
G.~E. Andrews and B.~C. Berndt, \emph{{R}amanujan's {L}ost {N}otebook, {P}art
  {I}}, Springer Science+Business Media, New York, 2005.

\bibitem{ramlet}
Bruce~C. Berndt and Robert~A. Rankin, \emph{Ramanujan: Letters and commentary},
  History of Mathematics, vol.~9, American Mathematical Society, Rhode Island,
  2000.

\bibitem{chokimkoo}
B.~Cho, N.~M. Kim, and J.~K. Koo, \emph{Affine models of the modular curves
  ${X}(p)$ and its application}, Ramanujan J. \textbf{24} (2011), 235--257.

\bibitem{chokoo}
B.~Cho and J.~K. Koo, \emph{Construction of class fields over imaginary
  quadratic fields and applications}, Q. J. Math. \textbf{61} (2010), 199--216.

\bibitem{chokoopark}
B.~Cho, J.~K. Koo, and Y.~K. Park, \emph{Arithmetic of the
  {R}amanujan-{G}\"{o}llnitz-{G}ordon continued fraction}, J. Number Theory
  \textbf{129} (2009), 922--947.

\bibitem{coop}
S.~Cooper, \emph{On {R}amanujan's function $k(q)=r(q)r(q^2)^2$}, Ramanujan J.
  \textbf{20} (2009), 311--328.

\bibitem{cooper}
S.~Cooper, \emph{Ramanujan's {T}heta {F}unctions}, Springer, Cham, 2017.

\bibitem{diashur}
F.~Diamond and J.~Shurman, \emph{A {F}irst {C}ourse in {M}odular {F}orms},
  Graduate Texts in Mathematics, vol. 228, Springer, New York, 2005.

\bibitem{gee}
A.~Gee, \emph{Class invariants by {S}himura's reciprocity law}, J. Th\'{e}or.
  Nombres Bordeaux \textbf{11} (1999), 45--72.

\bibitem{geehons}
A.~Gee and M.~Honsbeek, \emph{Singular values of the {R}ogers-{R}amanujan
  continued fraction}, Ramanujan J. \textbf{11} (2006), 267--284.

\bibitem{guad}
R.~Guadalupe, \emph{Modularity of a certain continued fraction of {R}amanujan},
  Ramanujan J. \textbf{63} (2024), 947--967.

\bibitem{isdishii}
N.~Ishida and N.~Ishii, \emph{The equations for modular function fields of
  principal congruence subgroups of prime level}, Manuscripta Math. \textbf{90}
  (1996), 271--285.

\bibitem{leepark}
Y.~Lee and Y.~K. Park, \emph{Ramanujan's function $k(\tau)=r(\tau)r(2\tau)^2$
  and its modularity}, Open Math. \textbf{18} (2020), 1727--1741.

\bibitem{ligo}
G.~Ligozat, \emph{Courbes modulaires de genus 1}, M\'{e}moires de la S. M. F.
  \textbf{43} (1975), 5--80.

\bibitem{newman1}
M.~Newman, \emph{Construction and application of a class of modular functions},
  Proc. London Math. Soc. \textbf{s3-7} (1957), 334--350.

\bibitem{newman2}
\bysame, \emph{Construction and application of a class of modular functions
  {(II)}}, Proc. London Math. Soc. \textbf{s3-9} (1959), 373--387.

\bibitem{ono}
K.~Ono, \emph{The {W}eb of {M}odularity: {A}rithmetic of the {C}oefficients of
  {M}odular {F}orms and $q$-series}, CBMS Regional Conference Series in
  Mathematics, vol. 102, American Mathematical Society, Rhode Island, 2004.

\bibitem{shim}
G.~Shimura, \emph{Introduction to the {A}rithmetic {T}heory of {A}utomorphic
  {F}unctions}, Princeton University Press, New Jersey, 1971.

\bibitem{wat}
G.~Watson, \emph{Theorems stated by {R}amanujan ({VII}): {T}heorems on
  continued fractions}, J. London Math. Soc. \textbf{4} (1929), 39--48.

\bibitem{xiayao}
E.~X.~W. Xia and O.~Y.~M. Yao, \emph{Some modular relations for {R}amanujan's
  function $k(q)=r(q)r^2(q^2)$}, Ramanujan J. \textbf{35} (2014), 243--251.

\bibitem{yang}
Y.~Yang, \emph{Transformation formulas for generalized {D}edekind eta
  functions}, Bull. Lond. Math. Soc. \textbf{36} (2004), 671--682.

\bibitem{yang2}
Y.~Yang, \emph{Defining equations of modular curves}, Adv. Math. \textbf{204}
  (2006), 481--508.

\bibitem{ye}
D.~Ye, \emph{Ramanujan's function $k$, revisited}, Ramanujan J. \textbf{56}
  (2021), 931--952.

\end{thebibliography}
\end{document}